\documentstyle[12pt]{article}
\begin{document}
\newcommand{\nc}{\newcommand}
\nc{\namelistlabel}[1]{\mbox{#1}\hfil}
\newenvironment{namelist}[1]{%
\begin{list}{}
{
\let\makelabel\namelistlabel
\settowidth{\labelwidth}{#1}
\setlength{\leftmargin}{1.1\labelwidth}
}
}{%
\end{list}}
\def\theequation{\thesection.\arabic{equation}}
\newtheorem{theorem}{\bf Theorem}[thetheorem]
\newtheorem{prop}{\bf Proposition}[theprop]
\newtheorem{corollary}{\bf Corollary}[thecorollary]
\newtheorem{remark}{\bf Remark}[theremark]
\newtheorem{lemma}{\bf Lemma}[thelemma]
\def\thetheorem{\thesection.\arabic{theorem}}
\def\theprop{\thesection.\arabic{prop}}
\def\theremark{\thesection.\arabic{remark}}
\def\thecorollary{\thesection.\arabic{corollary}}
\def\thelemma{\thesection.\arabic{lemma}}
\nc{\bsp}{\begin{sloppypar}}
\nc{\esp}{\end{sloppypar}}
\nc{\be}{\begin{equation}}
\nc{\ee}{\end{equation}}
\nc{\beanno}{\begin{eqnarray*}}
\nc{\inp}[2]{\left( {#1} ,\,{#2} \right)}
\nc{\dip}[2]{\left< {#1} ,\,{#2} \right>}
\nc{\disn}[1]{\|{#1}\|_h}
\nc{\pax}[1]{\frac{\partial{#1}}{\partial x}}
\nc{\tpar}[1]{\frac{\partial{#1}}{\partial t}}
\nc{\xpax}[2]{\frac{\partial^{#1}{#2}}{\partial x^{#1}}}
\nc{\pat}[2]{\frac{\partial^{#1}{#2}}{\partial t^{#1}}}
\nc{\ntpa}[2]{{\|\frac{\partial{#1}}{\partial t}\|}_{#2}}
\nc{\xpat}[2]{\frac{\partial^{#1}{#2}}{\partial t \partial x}}
\nc{\npat}[3]{{\|\frac{\partial^{#1}{#2}}{\partial t^{#1}}\|}_{#3}}
\nc{\xkpat}[3]{\frac{\partial^{#1}{#2}}{\partial t^{#3} \partial x}}
\nc{\jxpat}[3]{\frac{\partial^{#1}{#2}}{\partial t \partial x^{#3} }}
\nc{\eeanno}{\end{eqnarray*}}
\nc{\bea}{\begin{eqnarray}}
\nc{\eea}{\end{eqnarray}}
\nc{\ba}{\begin{array}}
\nc{\ea}{\end{array}}
\nc{\nno}{\nonumber}
\nc{\dou}{\partial}
\nc{\bc}{\begin{center}}
\nc{\ec}{\end{center}}
 \nc{\bb}{\mbox{\hspace{.25cm}}}
\nc{\bite}{\begin{itemize}}
\nc{\eite}{\end{itemize}}
\nc{\bth}{\begin{theorem}}
\nc{\eth}{\end{theorem}}
\nc{\bpr}{\begin{prop}}
\nc{\epr}{\end{prop}}
\nc{\blem}{\begin{lemma}}
\nc{\elem}{\end{lemma}}
\nc{\benu}{\begin{enumerate}}
\nc{\eenu}{\end{enumerate}}
\nc{\bcor}{\begin{corollary}}
\nc{\ecor}{\end{corollary}}
\nc{\Pf}{{\bf Proof. }}
\nc{\ddt}{\f{d}{dt}}
\nc{\ddr}{\f{d}{dr}}
\nc{\ddx}{\f{d}{dx}}
\nc{\DDx}{\f{d^{2}}{dx^{2}}}
\nc{\dodot}{\f{\dou}{\dou t}}
\nc{\dodori}{\f{\dou}{\dou r_i}}
\nc{\dodorj}{\f{\dou}{\dou r_j}}
\nc{\N}{I\!\!N}
\nc{\R}{I\!\!R}
\nc{\Hy}{I\!\!H}
\nc{\C}{I\!\!\!\!C}
\nc{\na}{\nabla}
\nc{\la}{\lambda}
\nc{\s}{\sinh}
\nc{\co}{\cosh}
\nc{\vl}{V_{\la}}
\nc{\ga}{\gamma}
\nc{\vg}{V_{\ga}}
\nc{\T}{\Theta}
\nc{\th}{\theta}
\nc{\f}{\frac}
\nc{\rw}{\rightarrow}
\nc{\lrw}{\Longrightarrow}
\nc{\om}{\omega}
\nc{\Om}{\Omega}
\nc{\al}{\alpha}
\nc{\ro}{\varrho}
\nc{\qed}{\phantom{a}{\hfill \rule{2.5mm}{2.5mm}}}
\nc{\D}{\Delta}
\nc{\TD}{\tilde{\D}}
\nc{\si}{\sigma}
\nc{\noi}{\noindent}
\nc{\bib}{\bibitem}
\title{\bf An Intrinsic Approach to Lichnerowicz Conjecture} 
\author{Akhil Ranjan }
\maketitle
\begin{abstract}
In this paper we give a proof of Lichnerowicz Conjecture for compact
simply connected manifolds which is intrinsic in the sense that it
avoids the {\it Nice Embeddings} into eigen spaces of the Laplacian.
Even if one wants to use these embeddings this paper gives a more
streamlined proof. 
\end{abstract}
{\bf Keywords:} Harmonic manifolds, Blaschke manifold, mean curvature, 
Jacobi differential equation, Ricci curvature, compact rank one symmetric
spaces, nice embeddings.
\footnotetext{AMS Subject Classification: 53C35}
\footnotetext{MRR 04-96}
\section{\bf Introduction}
The object of this paper is to present an intrinsic proof of the
Lichnerowicz's conjecture for the compact simply connected harmonic
manifolds. For the definition of harmonic manifolds see \cite{Bes1}. One of
characterisations is that the geodesic spheres around any point 
have constant mean curvature depending only on the radius of the
sphere. It suffices to consider small values of the radii. 
Lichnerowicz showed that for dimension less than or equal to 4, such
a manifold must be either flat or a locally symmetric space of
rank one (see \cite{Lic}, \cite{Bes1}). He quite naturally asked whether the same was 
true in higher dimensions. Great progress was made in the case of compact
simply connected harmonic manifolds, a detailed account of which is
given in Besse's book \cite{Bes1}. It is shown that these are all
Blaschke manifolds and their Ricci tensor is proportional to the
metric tensor or in other words they are Einstien manifolds.
From topological point of view each such manifold has its (integral)
cohomology ring isomorphic to precisely one of the compact rank one
symmetric space to be referred to as its {\it model CROSS} henceforth.
This result is due to Allemigeon.
A particularly striking discovery about compact harmonic spaces
is a family of isometric minimal immersions into the round spheres
in eigen-spaces of the Laplacian acting on the space of square-integrable 
functions. Moreover, any two geodesics were shown to be congruent to each other
under some Euclidean isometry. These are now known as {\it Besse's Nice
Embeddings}. In 1990 Szabo \cite{Sza1} successfully used them along with
other known facts about harmonic manifolds to answer Lichnerowicz's
query affirmatively for compact simply connected harmonic manifolds. 
In contrast to compact case, Damek and Ricci in 1992 \cite{DR1} 
produced a family of examples of homogeneous harmonic manifolds which are not 
locally symmetric. In  Szabo's paper the key point was to 
show that the volume function of a compact simply connected harmonic 
manifold when expressed in terms of normal coordinates coincided 
with that of its {\it model CROSS}. To this end he goes through the  
following steps:
\benu
\item He establishes what he calls {\it basic commutativity in harmonic
spaces}. This implies in conjunction with Allamigeon's theorem that
for any point $p$ on the manifold and any eigen value $\lambda$ of
the Laplacian there exist eigen-functions which depend only on
radial distance from $p$. Moreover, starting with any eigen function
and averageing over geodesic spheres around $p$ we get such a radial 
eigen function.
\item By moving the point $p$ along a geodesic and averageing in the
said manner we get a parallelly displaced family of functions in
the eigen-space which is finite dimensional. This along with the
fact that each geodesic is periodic of period {\it assumed} to be $2\pi$
enables one to conclude that the $radial$ eigen functions alluded
to above are polynomials in cosine of the radial distance.
\item At this stage the {\it nice embeddings} are used to pin down
the volume function in geodesic normal coordinates.
\item  Finally it follows that the first radial eigen function is
linear of the form $A\cos r + B$ and studying the nice embedding
in the first eigen-space one shows symmetry easily.
\eenu

In this paper we show that the {\it nice embeddings} can be avoided
in the step (3) above. Steps (1) and (2) do not require them 
anyway. Our proof of step (3) can be regarded as a streamlined
version of that given by Szabo. As for the last step one can either
use {\it nice embeddings} or the partial solution to a problem of
Antonio Ros about the first eigen value of {\em P}-manifolds given in 
\cite{RS}. In the latter case one has worked wholly within the manifold
thereby proving Szabo's theorem intrinsically.
\section{Laplacian on radial functions}
Let $\Theta (r)$ denote the volume function on a simply connected
compact harmonic manifold $M$ in terms of normal coordinates and
 $\sigma (r)$ the mean curvature of any geodesic sphere of radius $r$.
 It is easily shown that 
$$\frac{{\Theta}'(r)}{\Theta (r)} ~ =~ \sigma (r)$$
 Further for a point $p$ and an eigen-value $\lambda$ of the
Laplacian $\Delta$, let $u$ be an eigen function which depends only on
radial distance $r$ from $p$. As shown by Szabo (\cite{Sza1},p.8,eq.(2.1)) 
$u$ satifies the following ODE
\begin{equation}
\label{0de}
 u''+ {\sigma (r)}u'+ \lambda {u} = 0 
\end{equation}
Here $'$ means derivative with respect to $r$. We would like to
study how closely $\Theta $ and $\sigma $ agree with their anologues
in its {\it model CROSS}. Let us first define the volume function on all
of $real$ line as follows. 
 Consider a geodesic $\gamma$ through a point $p$. Let $J_2,...,J_d$ 
be the Jacobi fields along $\gamma$ which vanish at $\gamma (0)$ and whose 
derivatives at $\gamma (0)$ form an orthonormal basis along with 
$\gamma '(0)$. Let $E_1,...,E_d$ be parallel translation of the above
orthonormal basis along $\gamma$, $E_1$ being ${\gamma}'(r)$. Now set
$$\Theta (r) = <J_2\wedge ...\wedge J_d\:,\:E_2\wedge...\wedge E_d>(r)$$. 
By virtue of it being a Blaschke manifold, $\Theta$ when considered
as a function on whole of the real line has the following properties :
\benu
\item It is periodic of period $2\pi$.
\item It has zeroes of order $k-1$ at $r~=~n\pi $ for $n$ any odd integer
and zeroes of order $d-1$ at $r=n\pi $ for $n$ any even integer.
Here $d$ is the dimension of $M$ and $k$ is the degree of the 
generator of the cohomology ring of $M$.
\item  $\Theta (r)\:=\: (-1)^{d-1}\Theta (-r)$
\eenu
This clearly allows us to write
$$ \Theta (r)\; =\;e^{\alpha (\cos r)} \sin^{d-1}(r/2) \cos^{k-1}(r/2)$$
or $\Theta (r)\; =\;e^{\alpha (\cos r)}\Theta _0(r)$ where $\Theta _0(r)$
is the volume function of the {\it model CROSS} and $\alpha$ is a smooth
(actually analytic) function on [-1,1] with $\alpha (1)\; =\; 0$.
 Hence $\sigma (r)$ = ${\sigma_0(r)} -{\sin (r)\alpha '(\cos r)}$ {\rm since}
 $\sigma\;=\;\frac{\Theta '}{\Theta}$.\\
{\it Caution:} The conventional volume function is the absolute
value of the one we have defined. They both agree within the 
$injectivity~radius$ i.e. for $0 \leq r \leq \pi$. 
 For $0 < r < \pi$, an easy calculation gives that
$$ \sigma _0(r) \;=\;\f{1}{2\sin r}[(d-1)(1+\cos r)-(k-1)(1-\cos r)]$$
 By Lemma 4.2 of \cite{Sza1}, $u$ in eq.(2.1) is of the form $u\;=\;f(\cos r)$
for some {\it polynomial} $f$. Inserting all this data in $\ref{0de}$
and setting $\cos r = x$ we see that $f$ satisfies the following
\begin{equation} 
\label{0'de} 
(1-x^2)f''-[\f{d}{2}(1+x)-\f{k}{2}(1-x)+(1-x^2)\alpha'(x)]f'+\lambda f\;=\;0,
 -1\;\leq\;x\;\leq\;1.
\end{equation}
 In the above equation $'$ denotes derivative w.r.t. $x$.
\section{Jacobi Differential Equation}
 The differential equation
\be
\label{1de}
(1-x^2)u''-[(1+b)(1+x)-(1+a)(1-x)]u'+\la u\;=\;0
\ee
has been studied classically as a (singular) Sturm-Liouville equation
on [-1,1] and it is known that for $a$ and $b$ in $(-1,\infty)$ and
under natural boundary conditions ($u$ bounded as $|x| \rw 1$) solutions
exist for $\la = n(n+a+b+1), n \in \N$ and for each such value of
$\la$, $u$ is a polynomial of degree $n$. Moreover, $u$ is unique upto
a scalar multiple. In fact these polynomials form a complete orthogonal
system in $L^2([-1,1],\rho dx)$ where $\rho (x)\;=\;(1+x)^a(1-x)^b$ is the
weight function. These are known as Jacobi polynomials. (See \cite{BR} 
p. 289). This differential equation is known as Jacobi differential 
equation with parameters $a$ and $b$. We assume that $a , b > -1.$ 

 In this section we consider a {\it perturbed} Jacobi equation where
we have an extra term of the form $(1-x^2)\delta (x)$ as a coefficient
of $u'$, $\delta$ being a continuous function on [-1,1]. Comparing with
the corresponding equation satisfied by the polynomial $f$ in the 
previous section we easily see that
$1+b = \f{d}{2}$, $1+a = \f{k}{2}$ and $\delta = \alpha '$
 We also know that $k$ cannot exceed $d$ and can only take values in
${ 2,4,8,d}$, hence $a ,b > -1$ is clearly true. 
 Now we are ready to state our main theorem.
\begin{theorem}
 If the perturbed Jacobi differential equation 
\be
\label{2de} 
  (1-x^2)u''-[(1+b)(1+x)-(1+a)(1-x)+(1-x^2)\delta(x)]u'+\la u\; =\; 0
\ee
admits a nonconstant polynomial as a solution for some value of
 $\la$, then the perturbation term $\delta$ must vanish identically.
\end{theorem} 
\begin{corollary}
The perturbation term ${\alpha}'$ in $\ref{0'de}$ vainshes. Consequently
$\alpha$ is identically zero and hence $\sigma = \sigma _0$ as well as
$\Theta = \Theta _0$.
\end{corollary}
The proof of the above will be broken into two lemmas. Let $P$
be a polynomial which we assume to be nonconstant and monic which
satisfies $\ref{2de}$ for a suitable $\la$.
\begin{lemma}
$\delta$ must be a rational function with the degree of the numerator
being strictly less than that of the denominator.
\end{lemma}
{\bf Proof:} 
\beanno
\delta = \f{LP+\la P}{(1-x^2)P'}
\eeanno
where 
\beanno
L\ =\ (1-x^2)\DDx -[(1+b)(1+x)-(1+a)(1+x)]\ddx  
\eeanno
is the Jacobi differential operator (with parameters $a$ and $b$). Clearly
both numerator and denominator of $\delta$ are polynomials with
denominator nonvanishing and of degree strictly more than that of the
numerator. Hence the claim. 
\qed 
\begin{lemma}
Let $\delta = \f{p}{q}$ as a quotient of relatively prime polynomials
with q being monic. Then
\benu
\item All the roots of q are simple and in $\C \setminus [-1,1].$
\item $q|P'$ and $q|P$.
\item Let $q = \prod(x-\beta_i)$ and $m_i$ be the multiplicity of $\beta_i
\in P$, then $m_i \geq\,2$.
\item If we put $q_1 = \prod(x-\beta_i)^{m_i-1}$, then $\delta = \f{q_1'}{q_1}$.
\item $\pm$ 1 cannot be roots of P.
\item Roots of $P$ not common with those of $q$ are all simple.
\eenu
\end{lemma}
{\bf Proof:} Let $P = \prod(x-\beta_i)^{m_i}$ where $\beta_i$ are  
distinct complex numbers and $m_i$ are natural numbers which are
nonzero. Let 
\be
\label{v}
v = \f{P'}{P} = \sum \f{m_i}{x-\beta_i}
\ee
 Then $v$ satisfies the Riccati equation (see \cite{BR} p. 124) 
\be
\label{ric1}
 v'+v^2 = [\f{1+b}{1-x}-\f{1+a}{1+x}+\delta(x)]v-\f{\la}{1-x^2}
\ee
From \ref{v} we get
\be
\label{ric2}
v'+v^2 = \sum_i\f{m_i^2-m_i}{(x-\beta_i)^2} +\sum_{i\neq j}\f{2m_im_j}
{(\beta_i-\beta_j)(x-\beta_i)}
\ee
In the above equation we have expanded the {\it cross terms} occuring in
$v^2$ into partial fractions. Now let $q = \prod(x-\alpha_j)^{r_j}$ where
$\alpha_j$ are distinct complex numbers and $r_j \geq 1$. Since $p$ and
$q$ are relatively prime, if we expand $\delta = \f{p}{q}$ in partial
fractions, $\f{1}{(x-\alpha_j)^{r_j}}$ will survive for each $j$.
They will continue to survive after multiplication by $v = \sum \f{m_i}
{x-\beta_i}$ and further expansion into partial fractions. Now if we
compare the rhs of \ref{ric1} and \ref{ric2}  after expanding into
partial fractions we find that we must have $\{\alpha_j\} \subset \{\beta_i\}$
and $r_j = 1$ for each $j$. This proves that the roots of $q$ are simple.
Since $\delta =\f{p}{q}$ is continuous on [-1,1] roots of $q$ must 
be away from [-1,1]. This poves the first assertion.

$q|P$ is clear from above. To show that $q|P'$, we note that 
$LP+\la P = \f{p(1-x^2)P'}{q}$ is a polynomial. Hence $q|P'$ since it
is relatively prime to $p$, $1-x$, and $1+x$. This gives the second claim.

Put $S = \{\beta_i\}, S' = \{\alpha_j\}$, then $S' \subset S$. Also put 
$S'' = S \setminus S'$. 
We can then write $q = \prod_{S'}(x-\beta_j)$ and hence 
$\delta = \sum_{S'} \f{c_i}{x-\beta_i}$ where $c_i$ are nonzero numbers.
Coming back to the third 
statement, let us compare the coefficient of $\f{1}{(x-\beta_i)^2}$
in the rhs of \ref{ric1} and \ref{ric2} (after partial fractions)
for $\beta_i \neq \pm 1$. We see that 
\beanno
c_{i}\ m_{i}\ =\ m_{i}^{2}\ -\ m_{i}\ {\rm for}\ i\ {\rm s.t}\ \beta_{i}\
\in\ S'\ {\rm and}\ m_i^2-m_i= 0 \ \ {\rm if}\ \ \beta_i \in S''\setminus
\{\pm 1\} 
\eeanno
From this we can conclude that  
\beanno
c_{i}\ =\ m_{i}\ -\ 1 \ {\rm for}\ i\ \ s.t.\ \beta_i \in S'\ {\rm and}
\eeanno
\beanno
\ m_i\ =\ 1\ {\rm for}\ i\ s.t.\ \beta_i \in\ S''\setminus \{\pm 1\}\ ({\rm 
since}\  m_i\ \neq\ 0\ {\rm for each}\ i).
\eeanno
\beanno
\f{p}{q} = \sum_{S'}\f{m_i-1}{x-\beta_i} 
\eeanno 
The first of these shows that $m_i \geq 2$ for $\beta_i \in S'$ because
$c_i \neq 0$ and the second one can be rewritten as 
$\f{p}{q} = \f{q_1'}{q_1}$, where $q_1 = \prod_{S'}(x-\beta_i)^{m_i-1}$. 
These are just the third and the fourth assertions.
 
For the fifth claim, write $$P(x) = \prod_S(x-\beta_i)^{m_i} 
= (x-1)^A(x+1)^B\prod_{S'}(x-\beta_i)^{m_i}\prod_{S''\setminus 
\{\pm 1\}}(x-\beta_i) $$. 
Comparing coefficients of $(x-1)^{-2} and (x+1)^{-2}$ we see that
$$ A^2-A = -(1+b)A \  {\rm and}\  B^2-B = -(1+a)B$$
 Since $a$ and $b$ are more than -1 and $A$ and $B$ are natural
numbers we get $A = B =0$. Thus $S = S'\bigcup S''$ ; the set of roots 
of $P$, is disjoint from $\{\pm 1\}$

Finally, $S''\setminus \{\pm 1\}\ =\ S''$ so that for $\beta_i \in\ S''$
we have $m_i$ = 1 which is the sixth assertion.
\qed
\section{Proof of theorem 3.1} 
We have the following facts:
$P(x) = \prod_{S'}(x-\beta_i)^{m_i}\prod_{S''}(x-\beta_i), m_i \geq 2.
q_1(x) = \prod_{S'}(x-\beta_i)^{m_i-1}$.
Clearly, $q_1|P'$. Let $\f{P'}{q_1} = n\prod_j(x-\ga_j)^{M_j} = R(x)$,
say. (Here $n = degP$.) Then $\ga_j$ are those roots are $P'$ which
are not common with those of $P$. This is so because the roots of $P$
in $S''$ are all simple. Hence
\be
\label{eqn}
\{\ga_j\} \bigcap S = \emptyset = S\bigcap \{\pm 1\}.
\ee
Substituting the expressions obtained for $P, P'$ and $\delta$ in the
equation \ref{2de}, dividing by $(1-x^2)P'$ and simplifying we get
\be
\sum \f{M_j}{x-\ga_j} + \f{a+1}{x+1} + \f{b+1}{x+1} =
-\f{\la \prod_S(x-\beta_i)}{n\prod(x-\ga_j)^{M_j}}
\ee
 Putting $x = \beta \in S$ and using \ref{eqn} we find that
\be
\sum_j\f{M_j}{\beta -\ga_j} + \f{a+1}{\beta +1} +\f{b+1}{\beta -1}
= 0,\ \  {\rm for\ \  each \ \  such}\ \ \beta
\ee
This in turn implies that $\beta \in conv \{\ga_j, \pm 1 \}$ (where {\it conv}
denotes the {\it convex hull}). Or $S \subset conv \{\ga_j, \pm 1\}$. On
the other hand by Lucas' theorem $conv \{\ga_j \} \subset convS$. Now the
argument of Lemma 4.6 (\cite{Sza1}, p 23) goes through and shows that
$S \subset (-1,1)$. Also as observed earlier $S'$ is disjoint from
[-1,1]. Therefore, $S'$ must be empty so that $q$ is constant 1 and
hence $p =\delta$ vanishes identically as $deg(p) < deg(q)$. 
\qed
\section{Proof of the conjecture}
Let us first recall a theorem of Antonio Ros (\cite{Ros},Theorem 4.2, p.402)
\begin{theorem}
Let M be an n-dimensional $P_{2\pi}$-manifold and suppose that the Ricci
tensor, S, and the metric, <,>, on M verify the relation $S\geq k<,>,$
where k is a real constant. Let $\la_1$ be the first eigenvalue of the
Laplacian of M. Then we have
$$\la_1 \geq \f{1}{3} (2k+n+2)$$.
\end{theorem}
 He also remarked that for {\it CROSSES} the equality holds. He
naturally asked as to what restrictions apply to $M$ if equality
held.

As a partial answer to the above question the following theorem has
been proved (\cite{RS}) :
\begin{theorem}
If equality holds in the Ros' estimate for $\la_1$ of a P-manifold
and if M admits a corresponding eigen-function without saddle
points, then M is a CROSS.
\end{theorem}
Also see \cite{RS1} for another related result.
 
\noi{\bf Proof of the Conjecture:} 
Now from the corollary 3.1 to our main theorem $\Theta = \Theta_0$ and
hence $Ric_M = Ric_0$ and $\la_1(M) = \la_1(M_0)$ where $M_0$ 
denotes the {\it model CROSS}. Moreover, from any point on $M$
the first radial eigen-function is of the form $\cos r + C$ and
hence without saddle points. The proof of the Lichnerowicz 
conjecture is now complete.

\noi
Department of Mathematics\\
Indian Institute of Technology\\
Mumbai 400076, INDIA\\
email: aranjan@ganit.math.iitb.ernet.in 
\end{document}